\documentclass[11pt]{amsart}

\usepackage{latexsym}
\usepackage{amscd}
\usepackage[all]{xy}
\usepackage[mathscr]{euscript}
\usepackage{dsfont}

\normalsize

\RequirePackage{xspace}

\newcommand{\thought}[1]{}
\renewcommand{\thought}[1]{ \textbf{[#1]}}

\usepackage{enumerate}        
\newenvironment{roenumerate}{\begin{enumerate}[\upshape (i)]}{\end{enumerate}}

\newcommand\nc {\newcommand}
\newcommand\rnc{\renewcommand}

\newcount\blopone
\newcount\xone
\newcount\xtwo
\newcount\ytwo

\newtheorem{theorem}{Theorem}[section]
\newtheorem{prop}[theorem]{Proposition}
\newtheorem{com}[theorem]{Comment}
\newtheorem{apl}[theorem]{Application}
\newtheorem{exercise}[theorem]{Exercise}
\newtheorem{redu}[theorem]{Reduction}
\newtheorem{refinement}[theorem]{Refinement}
\newtheorem{summary}[theorem]{Summary}
\newtheorem{importnota}[theorem]{Important Notation}
\newtheorem{prblm}[theorem]{Problem}
\newtheorem{notation}[theorem]{Notation}
\newtheorem{explanation}[theorem]{Explanation}
\newtheorem{defin}[theorem]{Definition}
\newtheorem{caution}[theorem]{Caution}
\newtheorem{remark}[theorem]{Remark}
\newtheorem{reminder}[theorem]{Reminder}
\newtheorem{illustration}[theorem]{Illustration}
\newtheorem{observation}[theorem]{Observation}
\newtheorem{lemma}[theorem]{Lemma}
\newtheorem{construction}[theorem]{Construction}
\newtheorem{discussion}[theorem]{Discussion}
\newtheorem{corollary}[theorem]{Corollary}
\newtheorem{example}[theorem]{Example}
\newtheorem{conclusion}[theorem]{Conclusion}
\newtheorem{sketch}[theorem]{Sketch}
\newtheorem{triviality}[theorem]{Triviality}
\newtheorem{proto}[theorem]{Prototype Quasifibration}
\newtheorem{cauex}[theorem]{Cautionary Example}
\newtheorem{hypo}[theorem]{Hypothesis}
\newtheorem{subth}{ }[theorem]
\newtheorem{case}{Case}[theorem]
\newtheorem{ssubth}{ }[subth]
\newtheorem{facts}[theorem]{Facts}
\newtheorem{history}[theorem]{Historical Survey}
\newtheorem{proofs}[theorem]{Discussion of the Proofs, Old and New}

\nc\tri[1]{\begin{triviality}
\label{#1}}
\nc\fac[1]{\begin{facts}
\label{#1}
\begin{em}}
\nc\app[1]{\begin{apl}
\label{#1}
\begin{em}}
\nc\skt[1]{\begin{sketch}
\label{#1}
\begin{em}}
\nc\hst[1]{\begin{history}
\label{#1}
\begin{em}}
\nc\pfs[1]{\begin{proofs}
\label{#1}
\begin{em}}
\nc\cas[1]{\begin{case}
\label{#1}
\begin{em}}
\nc\rfn[1]{\begin{refinement}
\label{#1}}
\nc\prt[1]{\begin{proto}
\label{#1}}
\nc\lem[1]{\begin{lemma}
\label{#1}}
\nc\pro[1]{\begin{prop}
\label{#1}}
\nc\thm[1]{\begin{theorem}
\label{#1}}
\nc\dis[1]{\begin{discussion}
\label{#1}
\begin{em}}
\nc\cor[1]{\begin{corollary}
\label{#1}}
\nc\dfn[1]{\begin{defin}
\label{#1}}
\nc\sthm[1]{\begin{subth}
\label{#1}}
\nc\exm[1]{\begin{example}
\label{#1}
\begin{em}}
\nc\obs[1]{\begin{observation}
\label{#1}
\begin{em}}
\nc\plm[1]{\begin{prblm}
\label{#1}
\begin{em}}
\nc\rmk[1]{\begin{remark}
\label{#1}
\begin{em}}
\nc\rmd[1]{\begin{reminder}
\label{#1}
\begin{em}}
\nc\ntn[1]{\begin{notation}
\label{#1}
\begin{em}}
\nc\exe[1]{\begin{exercise}
\label{#1}
\begin{em}}
\nc\xpl[1]{\begin{explanation}
\label{#1}
\begin{em}}
\nc\smr[1]{\begin{summary}
\label{#1}
\begin{em}}
\nc\cau[1]{\begin{caution}
\label{#1}
\begin{em}}
\nc\hyp[1]{\begin{hypo}
\label{#1}
\begin{em}}
\nc\imn[1]{\begin{importnota}
\label{#1}
\begin{em}}
\nc\rdn[1]{\begin{redu}
\label{#1}
\begin{em}}
\nc\cax[1]{\begin{cauex}
\label{#1}
\begin{em}}
\nc\cmt[1]{\begin{com}
\label{#1}
\begin{em}}
\nc\con[1]{\begin{construction}
\label{#1}
\begin{em}}
\nc\ill[1]{\begin{illustration}
\label{#1}
\begin{em}}
\nc\ssthm[1]{\begin{ssubth}
\label{#1}
\begin{em}}
\nc\cnc[1]{\begin{conclusion}
\label{#1}
\begin{em}}

\nc\elem{\end{lemma}}
\nc\erdn{\end{em}\end{redu}}
\nc\erfn{\end{refinement}}
\nc\eprt{\end{proto}}
\nc\ethm{\end{theorem}}
\nc\ecor{\end{corollary}}
\nc\edfn{\end{defin}}
\nc\esthm{\end{subth}}
\nc\epro{\end{prop}}
\nc\etri{\end{triviality}}
\nc\eexm{\end{em}
\end{example}}
\nc\eobs{\end{em}
\end{observation}}
\nc\ecmt{\end{em}
\end{com}}
\nc\efac{\end{em}
\end{facts}}
\nc\eapp{\end{em}
\end{apl}}
\nc\ermk{\end{em}
\end{remark}}
\nc\ermd{\end{em}
\end{reminder}}
\nc\eill{\end{em}
\end{illustration}}
\nc\eplm{\end{em}
\end{prblm}}
\nc\ecas{\end{em}
\end{case}}
\nc\eskt{\end{em}
\end{sketch}}
\nc\ecau{\end{em}
\end{caution}}
\nc\ecax{\end{em}
\end{cauex}}
\nc\eimn{\end{em}
\end{importnota}}
\nc\entn{\end{em}
\end{notation}}
\nc\eexe{\end{em}
\end{exercise}}
\nc\expl{\end{em}
\end{explanation}}
\nc\edis{\end{em}
\end{discussion}}
\nc\econ{\end{em}
\end{construction}}
\nc\esmr{\end{em}
\end{summary}}
\nc\ehst{\end{em}
\end{history}}
\nc\epfs{\end{em}
\end{proofs}}
\nc\ehyp{
\end{em}
\end{hypo}}
\nc\ecnc{\end{em}
\end{conclusion}}
\nc\essthm{\end{em}
\end{ssubth}}

\nc\sst{\scriptstyle}
\newcommand{\comment}[1]{}
\newcommand{\ri}{\longrightarrow}

\newcommand{\zz}{{\mathbb Z}}

\newcommand{\nn}{{\mathbb N}}

\newcommand{\K}{{\mathbf K}}
\newcommand{\A}{\mathbf{Ac}}
\newcommand{\D}{{\mathbf D}}

\newcommand{\qq}{{\mathbb Q}}

\nc\op{^{\hbox{\rm\tiny op}}}
\nc\mth{^{\hbox{\rm\tiny th}}}

\nc\script{\mathscr}
\nc\z{\zeta}
\nc\bc{{\mathbb{BC}}}
\nc\KK{{\mathbb{K}}}
\nc\ct{{\script T}}
\nc\cf{{\script F}}
\nc\cg{{\script G}}
\nc\ch{{\script H}}
\nc\ck{{\script K}}
\nc\cl{{\script L}}
\nc\cm{{\script M}}
\nc\cv{{\script V}}
\nc\ce{{\script E}}
\nc\cs{{\script S}}
\nc\car{{\script R}}
\nc\cd{{\script D}}
\nc\cc{{\script C}}
\nc\ca{{\script A}}
\nc\ci{{\script I}}
\nc\cj{{\script J}}
\nc\co{{\script O}}
\nc\cu{{\script U}}
\nc\cw{{\script W}}
\nc\cx{{\script X}}
\nc\Cp{{\script P}}
\nc\cq{{\script Q}}
\nc\cy{{\script Y}}
\nc\cz{{\script Z}}
\nc\bd{\begin{description}}
\nc\ed{\end{description}}
\nc\ctob{{\script C}at\big(\ci^{op},\ca\big)}
\nc\clim{{\ds\mathop{\rm lim}_{\ds\longleftarrow}}\,}
\nc\climi{\clim_{\!i}\,}
\nc\climn{\clim^{\!n}\,}
\nc\colim{{\ds\mathop{\rm colim}_{\ds\la}}}
\nc\colimj{{\ds\mathop{\rm colim}_{\ds\la}}{}_{j\,}}
\nc\oa{\overline{\ca}}
\nc\s{\sigma}
\nc\ta{\tau}
\nc\os{\overline\sigma}
\nc\ot{\overline\tau}
\nc\T{\Sigma}
\nc\Tm{\Sigma^{-1}}
\nc\de[1]{{\mathop{\rm deg(#1)}}}
\nc\Ad[1]{\mathop{\rm Ad}(#1)}
\nc\ad[1]{\mathop{\rm ad}(#1)}
\nc\kth{{\it K}--theory}
\nc\loc[1]{{\text{\rm Loc}(#1)}}
\nc\coloc[1]{{\text{\rm Coloc}(#1)}}

\def\der #1 {D\left(#1\right)}
\nc\prf{\begin{proof}}
\nc\eprf{\end{proof}}
\nc\ds{\displaystyle}
\nc\Tor{\text{\rm Tor}}

\nc\cb{{\script B}}
\nc\ab{{\script A}b}

\nc\be{\begin{roenumerate}}
\nc\ee{\end{roenumerate}}

\nc\cat[1]{{\script C}at\Big({\big\{#1\big\}}\op\,\,,\,\,\ab\Big)}
\nc\csab{{\script C}at\big(\cs^{op},\ab\big)}
\nc\ctab{{\script C}at\Big({\{\ct^\alpha\}}^{op},\ab\Big)}
\nc\csex{{\script E}x\big(\cs^{op},\ab\big)}
\nc\ctex{{\script E}x\Big({\{\ct^\alpha\}}^{op},\ab\Big)}
\nc\sub{\qquad\subset\qquad}
\nc\ctr[1]{{\left.\ct\left(-,#1\right)\right|}_{\cs}}
\nc\ctrf[2]{{\left.\ct\left(#1,#2\right)\right|}_{\cs}}
\nc\Ctr[1]{{\left.\ct\left(-,#1\right)\right|}_{\ct^\alpha}}
\nc\Ctrf[2]{{\left.\ct\left(#1,#2\right)\right|}_{\ct^\alpha}}

\nc\la{\longrightarrow}
\nc\nin{\noindent}
\nc\cad[1]{\text{card}(#1)}
\nc\eq{\quad=\quad}
\nc\BA{\begin{array}{c}}
\nc\EA{\end{array}}
\nc\barr{
\[
\begin{array}{cccccccccccccccc}
}
\nc\earr{
\end{array}
\]
}
\nc\as[1]{{\langle S\rangle}^{#1}}
\nc\sh{\text{\it shift}}

\nc\yy[1]{{\left.\ct\left(-,#1\right)\right|}_{\ct^c}}
\nc\vrep[2]{{\left.\ct\left(#1,#2\right)\right|}_{\ct^\alpha}}
\nc\da{\downarrow}
\nc\Hom{{\mathop{\rm Hom}}}
\nc\HHom{{\script H}{\mathop{\rm om}}}
\nc\End{{\mathop{\rm End}}}
\nc\Ext{{\mathop{\rm Ext}}}
\nc\mr\Modtc
\nc\PExt{{\mathop{\rm PExt}}}
\nc\stm{\text{\rm stmod}(kG)}
\nc\stM{\text{\rm StMod}(kG)}
\nc\e{\varepsilon}
\nc\p{\varphi}

\nc\rs{\s^{-1}A}
\nc\br{{\{\s^{-1}A\}}}
\nc\ra\ri
\nc\y[1]{\mathbf{y}#1}
\nc\x[1]{\mathbf{z}#1}
\nc\mmod[1]{\text{\rm mod--}#1}
\nc\Mod[1]{#1\text{--\rm Mod}}
\nc\Md {\ensuremath{\mathop{\textup{Mod}}}}
\rnc\mod[1]{\ensuremath{\mathop{#1\textup{--mod}}}\xspace}
\nc\MMod[1]{\text{Mod-}#1}
\nc\Modtc{\Mod{\ct^c}}
\nc\pgldim[1]{\mathop{\rm pgldim}\,#1}
\nc\tf{{\rm [TR5]}}
\nc\tfs{{\rm [TR5$^*$]}}
\nc\Fun{\text{\rm Funct}(F\op,\ab)}
\nc\sym{\text{\rm Sym}}
\nc\sgn{\text{\rm sgn}}
\nc\Pro{\text{\rm Prod}^{}_\alpha(F\op,\ab)}
\nc\Yt[1]{{\left.\Hom_\ct^{}\left(-,#1\right)\right|}_F^{}}
\nc\dl{\delta}
\nc\Proj[1]{#1\text{--\rm Proj}}
\nc\proj[1]{#1\text{--\rm proj}}
\nc\Flat[1]{#1\text{--\rm Flat}}
\nc\Inj[1]{#1\text{--\rm Inj}}
\nc\Ima{\mathrm{Im}}
\nc\Ker{\mathrm{Ker}}
\nc\ov{\overline}
\nc\wt{\widetilde}
\nc\wh{\widehat}
\nc\ph{\varphi}
\nc\tstr{{\it t}--structure}
\nc\tstrs{{\it t}--structures}
\nc\spec[1]{{\text{\rm Spec}(#1)}}

\nc\EProd{\text{\rm EProd}}
\nc\ECoprod{\text{\rm ECoprod}}
\nc\Prod{\text{\rm Prod}}
\nc\ldimp{\text{\rm LDim}^{\prod}}
\nc\ldimc{\text{\rm LDim}^{\coprod}}
\nc\gen[2]{{\langle#1\rangle}^{}_{#2}}
\nc\Gen[2]{{\big\langle#1\big\rangle}^{}_{#2}}
\nc\genu[3]{{\langle#1\rangle}^{[#3]}_{#2}}
\nc\ogen[1]{\ov{\langle#1\rangle}}
\nc\ogenun[2]{\ov{\langle#1\rangle}_{#2}^{}}
\nc\ogenu[3]{\ov{\langle#1\rangle}^{[#3]}_{#2}}
\nc\ogenul[3]{\ov{\langle#1\rangle}^{(-\infty,#3]}_{#2}}
\nc\ogenuf[3]{\ov{\langle#1\rangle}^{[#3,\infty)}_{#2}}
\nc\genuf[3]{{\langle#1\rangle}^{[#3,\infty)}_{#2}}
\nc\genul[3]{{\langle#1\rangle}^{(-\infty,#3]}_{#2}}
\nc\dperf[1]{\D^{\mathrm{perf}}(#1)}
\nc\dcoh{\mathbf{D}^b_{\mathrm{coh}}}

\nc\dmcoh{\mathbf{D}^-_{\mathrm{coh}}}
\nc\dscoh{\mathbf{D}^{}_{\mathrm{coh}}}
\nc\RHHom{{\script{RH}}{\mathrm{om}}}
\nc\Coprod{\mathrm{Coprod}}
\nc\COprod{\mathrm{coprod}}
\nc\add{\mathrm{add}}
\nc\Add{\mathrm{Add}}
\nc\Smr{\mathrm{smd}}
\nc\id{\mathrm{id}}
\nc\LL{\mathbf{L}}
\nc\R{\mathbf{R}}
\nc\CC{\mathbf{C}}
\nc\wi{\wt{\text{\it\i}}}
\nc\exal{\ce\text{\it x}_\alpha(\ct^\alpha,\ab)}
\nc\exalz{\ce\text{\it x}_{\aleph_0}^{}(\ct^\alpha,\ab)}

\nc\tst[1]{\left({#1}^{\leq0},{#1}^{\geq0}\right)}
\nc\perf[1]{\script{P}\mathrm{erf}\left(#1\right)}
\nc\SF[1]{\mathrm{SF}\left(#1\right)}

\nc\one{\mathds{1}}
\nc\fc{\mathfrak{C}}
\nc\fl{\mathfrak{L}}
\nc\fs{\mathfrak{S}}
\nc\Prf{\text{\bf Perf}}
\nc\qc{\text{\bf qc}}
\nc\vect[1]{\script{V}\mathit{ect}\left(#1\right)}
\nc\vectgr[1]{\script{V}\mathit{ect}\script{G}\mathit{r}\left(#1\right)}
\nc\Ch{\mathbf{Ch}}

\nc\hoco{
\begin{picture}(40,10)
\put(20,0){\makebox(0,0)[b]{\text{\rm Hocolim}}}
\put(5,-2){\vector(1,0){30}}
\end{picture}\,\,}

\nc\holim{{
\begin{picture}(40,10)
\put(20,0){\makebox(0,0)[b]{\text{\rm Holim}}}
\put(35,-2){\vector(-1,0){30}}
\end{picture}}}

\begin{document}

\author{Amnon Neeman}\thanks{The research was partly supported 
by the Australian Research Council}
\address{Centre for Mathematics and its Applications \\
        Mathematical Sciences Institute\\
        Building 145\\
        The Australian National University\\
        Canberra, ACT 2601\\
        AUSTRALIA}
\email{Amnon.Neeman@anu.edu.au}

\title{A counterexample to vanishing conjectures for negative $K$-theory}

\begin{abstract}
In a 2006 article Schlichting conjectured that the negative {\it K--}theory
of any abelian category must vanish. This conjecture was generalized
in a 2019 article by Antieau, Gepner and Heller, who hypothesized 
that the negative  {\it K--}theory of any category with a bounded
{\it t--}structure must vanish.

Both conjectures will be shown to be false.
\end{abstract}



\maketitle

\tableofcontents

\setcounter{section}{-1}

\section{Introduction}
\label{S0}

Let $\ce$ be any idempotent-complete exact category. 
We may form the category $\A^b(\ce)$,
whose objects are the acyclic bounded cochain complexes
in $\ce$. In Lemma~\ref{L7.3} and Remark~\ref{R47.9957}
we prove that $\A^b(\ce)$ always
has a bounded \tstr. We eventually show that, if
$\ce=\vect Y$ is the category of vector bundles on a projective
curve $Y$ with only simple nodes as singularities, then 
there is an injective map $\KK_{-1}(\ce)\la\KK_{-2}\big[\A^b(\ce)\big]$.
Since there are known examples of nodal curves for which
$\KK_{-1}\big[\vect Y\big]\neq0$, this provides
a counterexample to 
Antieau, Gepner and Heller's~\cite[Conjecture~B]{Antieau-Gepner-Heller19}.

More generally: let $\ce$ be any idempotent-complete exact
category. In Proposition~\ref{P8.1} we produce a
homotopy fiber sequence
\[\xymatrix@C+20pt{
\KK\big[\A^b(\ce)\big]\ar[r] & \KK(\ce^\oplus) \ar[r] & \KK(\ce)
}\]
where $\ce^\oplus$ is $\ce$ with the split exact structure. Thus
the vanishing of $\KK_n\big[\A^b(\ce)\big]$ for
$n<0$ would imply that the natural maps 
$\KK_{n}(\ce^\oplus) \la \KK_n(\ce)$ must be isomorphisms 
for all $n<0$. It is entirely possible that there are 
many more counterexamples out there; the one computed in this 
article is the case of projective nodal curves. More
precisely: if $\ce=\vect Y$ with $Y$ a projective nodal curve,
we prove that $\KK_{-1}(\ce^\oplus)=0$. But there are 
known examples where $\KK_{-1}(\ce)\neq0$.

Following 
Beilinson, Bernstein and Deligne~\cite[Proposition~3.1.10]{BeiBerDel82},
if a triangulated category $\ct$ with a \tstr\ comes from a model and
has suitable ``filtered'' versions,
then there is a natural functor 
$F:\D^b\big[\ct^\heartsuit\big]\la\ct^b$, from the 
bounded derived
category of the heart of $\ct$ to the bounded
part $\ct^b\subset\ct$. And what's important
here is that the proof goes by a way that lifts to models.
If we apply this to $\A^b(\ce)$ we deduce
an induced map in \kth\ of the form
$\KK\big[\A^b(\ce)^\heartsuit\big]\la\KK\big[\A^b(\ce)\big]$.
And to show that the map in \kth\ is a homotopy equivalence, it suffices
to prove that 
$F:\D^b\big[\A^b(\ce)^\heartsuit\big]\la\A^b(\ce)$
is an equivalence of triangulated categories.
In Proposition~\ref{P47.3} we show that the 
functor 
$F:\D^b\big[\A^b(\ce)^\heartsuit\big]\la\A^b(\ce)$ 
is an equivalence if and only if the exact category $\ce$
is hereditary, meaning $\Ext^i(E,E')=0$ for all $E,E'\in\ce$
and $i\geq2$.

Since the category of vector bundles on a curve is hereditary,
we deduce that in our counterexample the 
map $F$ is an equivalence. Therefore
$\KK_{-2}\big[\A^b(\ce)^\heartsuit\big]\neq0$, giving
a counterexample to 
Schlichting~\cite[Conjecture~1 of Section~10]{Schlichting06},
which is also
Antieau, Gepner and Heller~\cite[Conjecture~A]{Antieau-Gepner-Heller19}.

\medskip

\nin
{\bf Acknowledgements.}\ \ The author would like to thank
Ben Antieau, Igor Burban, Ching-Li Chai, Xiao-Wu Chen,
Bernhard Keller, Henning Krause, Peter Newstead, 
Sundararaman Ramanan and Chuck Weibel 
for helpful comments and improvements
on earlier incarnations of the manuscript.

\section{The \tstr\ on the category $\A^b(\ce)$}
\label{S7}

\ntn{N7.1}
Let $\ce$ be an
idempotent-complete exact category, and let $\K(\ce)$ be the category
whose objects are the cochain
complexes of objects in $\ce$ and whose 
morphisms are the homotopy equivalence classes of cochain maps. Let
$\A(\ce)$ be the full subcategory of acyclic 
complexes. The full subcategories $\A^-(\ce)\subset\K^-(\ce)$,
$\A^+(\ce)\subset\K^+(\ce)$ and
$\A^b(\ce)\subset\K^b(\ce)$
are the obvious bounded versions\footnote{In the Introduction we followed
the notation of Schlichting~\cite{Schlichting06}, where
$\A^b(\ce)$ is a model category. In almost all of
the article we will follow the notation of
Krause~\cite{Krause15}, where $\A^b(\ce)$ stands for
the associated triangulated category.}. We remind the
reader of the definition of acyclicity: a cochain complex 
\[\xymatrix@C+20pt{
\cdots \ar[r]^-{\partial^{i-2}} & E^{i-1}\ar[r]^-{\partial^{i-1}} &
E^{i}\ar[r]^-{\partial^{i}} &
E^{i+1}\ar[r]^-{\partial^{i+1}} &\cdots
}\] 
possibly bounded, is declared
\emph{acyclic} if there exist admissible short exact sequences
\[\xymatrix@C+20pt{
0 \ar[r] & K^{i}\ar[r]^-{\alpha^{i}} &
E^{i}\ar[r]^-{\beta^{i}} &
K^{i+1}\ar[r] &0
}\] 
such that $\partial^i=\alpha^{i+1}\circ\beta^i$. 
The derived categories $\D^?(\ce)$
are defined to be the Verdier quotients $\K^?(\ce)/\A^?(\ce)$,
for $?$ being $b$, $-$, $+$ or the empty restriction.

Note that we are assuming $\ce$ idempotent-complete, and
\cite[Lemma~1.2]{Neeman90} proves that $\A(\ce)$ is a thick
subcategory of $\K(\ce)$. The fact that
$\A^-(\ce)\subset\K^-(\ce)$,
$\A^+(\ce)\subset\K^+(\ce)$ and
$\A^b(\ce)\subset\K^b(\ce)$
are all thick subcategories is older, it may essentially
be found in
Thomason and Trobaugh~\cite[1.11.1 (see also Appendix A)]{ThomTro}.
See also
\cite[Remark~1.10]{Neeman90} for a brief synopsis of 
the argument in Thomason-Trobaugh.

We will usually write $E^*$ as a shorthand for the object
\[\xymatrix@C+20pt{
\cdots \ar[r]^-{\partial^{i-2}} & E^{i-1}\ar[r]^-{\partial^{i-1}} &
E^{i}\ar[r]^-{\partial^{i}} &
E^{i+1}\ar[r]^-{\partial^{i+1}} &\cdots
}\]
in $\K(\ce)$.
\entn

\lem{L7.3}
Let $\ce$ be an idempotent-complete exact category, and 
let $\A(\ce)$ be the subcategory of
acyclics as in Notation~\ref{N7.1}. Define the full subcategories
\begin{eqnarray*}
\A(\ce)^{\leq0} &=& \{E^*\in\A(\ce)\mid E^i=0\text{ for all }i>0\} \\
\A(\ce)^{\geq0} &=& \{E^*\in\A(\ce)\mid E^i=0\text{ for all }i<-2\}
\end{eqnarray*}
Then the pair $\big[\A(\ce)^{\leq0}\,,\,\A(\ce)^{\geq0}\big]$ define
a \tstr\ on $\A(\ce)$.
\elem

\prf
The containments $\T\A(\ce)^{\leq0}\subset\A(\ce)^{\leq0}$
and $\A(\ce)^{\geq0}\subset\T\A(\ce)^{\geq0}$ are obvious from
the definition.

Now suppose we are given a morphism from an object $E^*\in\A(\ce)^{\leq0}$
to an object $F^*\in\A(\ce)^{\geq1}$. We may represent it by a cochain map
\[\xymatrix@C+20pt{
\cdots \ar[r]^-{\partial^{-3}} & E^{-2}\ar[r]^-{\partial^{-2}}\ar[d] &
E^{-1}\ar[r]^-{\partial^{-1}}\ar[d]_{f} &
E^{0}\ar[r]\ar[d]^g &0\ar[r]\ar[d] &\cdots\\
\cdots \ar[r] & 0\ar[r] &
F^{-1}\ar[r]^-{\wt\partial^{-1}} &
F^{0}\ar[r]^-{\wt\partial^{0}} & F^1\ar[r]^-{\wt\partial^{1}} &
\cdots
}\]
The fact that $\wt\partial^0\circ g=0$ says that $g$ must factor 
uniquely through
the kernel of the map $\wt\partial^{0}$, which happens to be
the map $\wt\partial^{-1}:F^{-1}\la F^0$. Thus we may find a (unique)
morphism $\theta:E^0\la F^{-1}$ with $g=\wt\partial^{-1}\circ\theta$.
But now we have the equalities
\begin{eqnarray*}
\wt\partial^{-1}\circ f &=& g\circ\partial^{-1} \\
&=& \wt\partial^{-1}\circ\theta\circ\partial^{-1}
\end{eqnarray*}
where the first comes from the commutativity implied by the cochain map
$E^*\la F^*$,
and the second is by precomposing $g=\wt\partial^{-1}\circ\theta$ with
$\partial^{-1}$.
And, since $\wt\partial^{-1}$ is a monomorphism (even an admissible
monomorphism), it follows that $f=\theta\circ\partial^{-1}$. Thus $\theta$
provides a homotopy of the cochain map $E^*\la F^*$ with the zero map.

Next choose any object $E^*\in\A(\ce)$, that is a complex
\[\xymatrix@C+20pt{
\cdots \ar[r]^-{\partial^{i-2}} & E^{i-1}\ar[r]^-{\partial^{i-1}} &
E^{i}\ar[r]^-{\partial^{i}} &
E^{i+1}\ar[r]^-{\partial^{i+1}} &\cdots
}\]
such that each morphism $\partial^i:E^i\la E^{i+1}$ has a factorization 
$E^i\stackrel{\beta^i}\la K^{i+1}\stackrel{\alpha^{i+1}}\la E^{i+1}$ as
in Notation~\ref{N7.1}. In particular: we may write 
$\partial^{-1}:E^{-1}\la E^0$ as a composite 
$E^{-1}\stackrel{\beta^{-1}}\la K^{0}\stackrel{\alpha^{0}}\la E^{0}$. But 
now consider the cochain maps
\[\xymatrix@C+20pt{
\cdots \ar[r]^-{\partial^{-3}} & E^{-2}\ar[r]^-{\partial^{-2}}\ar[d]^\id &
E^{-1}\ar[r]^-{\beta^{-1}}\ar[d]^\id &
K^{0}\ar[r]\ar[d]^{\alpha^0} &
0\ar[r]\ar[d] &\cdots\\
\cdots \ar[r]^-{\partial^{-3}} & E^{-2}\ar[r]^-{\partial^{-2}}\ar[d] &
E^{-1}\ar[r]^-{\partial^{-1}}\ar[d]^{\beta^{-1}} &
E^{0}\ar[r]^-{\partial^{0}}\ar[d]^\id &
E^{1}\ar[r]^-{\partial^{1}}\ar[d]^\id &\cdots\\
\cdots \ar[r] & 0\ar[r]\ar[d] &
K^{0}\ar[r]^-{\alpha^{0}}\ar[d]^\id &
E^{0}\ar[r]^-{\partial^{0}}\ar[d] &
E^{1}\ar[r]^-{\partial^{1}}\ar[d] &\cdots\\
\cdots \ar[r]^-{-\partial^{-2}} &
E^{-1}\ar[r]^-{-\beta^{-1}} &
K^{0}\ar[r] &
0\ar[r] &0\ar[r] &\cdots\\
}\]
and we leave it to the reader to check that this is isomorphic
in $\A(\ce)$ to a distinguished
triangle $A^*\la E^*\la B^*\la\T A^*$, in 
which obviously $A^*\in\A(\ce)^{\leq0}$ and
$B^*\in\A(\ce)^{\geq1}$.

This completes the proof that the pair 
$\big[\A(\ce)^{\leq0}\,,\,\A(\ce)^{\geq0}\big]$ define
a \tstr\ on $\A(\ce)$.
\eprf

\rmk{R47.9957}
Given a triangulated category $\ct$ with a \tstr, it is customary to define
the subcategories
\[
\ct^-=\bigcup_{n=1}^\infty\ct^{\leq n}\,,\qquad\qquad
\ct^-=\bigcup_{n=1}^\infty\ct^{\geq -n}\,,\qquad\qquad
\ct^b=\ct^-\cap\ct^+\,.
\]
In the special case where $\ct=\A(\ce)$ and the \tstr\ is
as in Lemma~\ref{L7.3}, the definitions give that
\[
\left[\A(\ce)\right]^-=\A^-(\ce)\,,\qquad\qquad
\left[\A(\ce)\right]^+=\A^+(\ce)\,,\qquad\qquad
\left[\A(\ce)\right]^b=\A^b(\ce)\,,
\]
with $\A^?(\ce)$ being as in Notation~\ref{N7.1}. It follows that
the \tstr\ on $\A(\ce)$ restricts to {\it t---}structures on
$\A^?(\ce)$, with $?$ being each of $-$, $+$ and $b$.
And all four categories have the same heart.
\ermk

\rmk{R7.5}
The heart of the \tstr\ of Lemma~\ref{L7.3} is, by definition, given 
by the formula
$\A(\ce)^\heartsuit=\A(\ce)^{\leq0}\cap\A(\ce)^{\geq0}$, and
the formula gives that the objects  of $\A(\ce)^\heartsuit$
are the admissible short exact sequences
\[\xymatrix@C+20pt{
0\ar[r] & E^{-2}\ar[r] & E^{-1}\ar[r] & E^0 \ar[r] & 0
}\]
in the category $\ce$. Since $\A(\ce)^\heartsuit$ is a full subcategory
of $\K(\ce)$, the morphisms in $\A(\ce)^\heartsuit$
are homotopy equivalence classes
of cochain maps 
\[\xymatrix@C+20pt{
0\ar[r] & E^{-2}\ar[r]\ar[d] & E^{-1}\ar[r]\ar[d] & E^0 \ar[r]\ar[d] & 0\\
0\ar[r] & F^{-2}\ar[r] & F^{-1}\ar[r] & F^0 \ar[r] & 0
}\]
With this description it isn't immediately obvious where this construction
comes from, let alone why this category must be abelian.
\ermk

\rmk{R7.6}
The abelian category $\A(\ce)^\heartsuit$ isn't new, it may be found in
Schlichting~\cite[Lemma~9 of Section~11]{Schlichting06}. In Schlichting's
presentation this category doesn't come as the heart of some 
\tstr, instead it is described as a subcategory of the category
$\text{Eff}(\ce)\subset\MMod\ce$, whose objects
are the effaceable functors in
the category $\MMod\ce=\Hom\big(\ce\op,\ab\big)$ 
of additive functors $\ce\op\la\ab$.
\ermk

\rmk{R7.7}
It might help to consider the special case 
where $\ce$ is an abelian category.
The Yoneda map $Y:\ce\la\mmod\ce$ embeds $\ce$ fully faithfully into the 
category $\mmod\ce\subset\Hom\big(\ce\op,\ab\big)$ 
of finitely presented functors $\ce\op\la\ab$. Recall: a functor 
$F:\ce\op\la\ab$ is \emph{finitely presented} if 
there exists an exact sequence
\[\xymatrix@C+20pt{
Y(A)\ar[r]^-{Y(f)} & Y(B) \ar[r] & F\ar[r] & 0
}\]
which can be thought of as a finite presentation of $F$ 
in the abelian category $\Hom\big(\ce\op,\ab\big)$.
Auslander's work  tells us
that the functor $Y:\ce\la\mmod\ce$ has an exact left adjoint 
$\Lambda:\mmod\ce\la\ce$. The way to compute 
$\Lambda(F)$ is to choose a finite 
presentation as above, and define 
$\Lambda(F)$ to be the cokernel of the map
$f:A\la B$.
With $\text{eff}(\ce)$ defined to be 
the full subcategory of $\mmod\ce$ annihilated
by the functor $\Lambda$, Auslander's 
formula~\cite[page~205]{Auslander66} 
goes on to tell us that
$\ce$ is the Gabriel quotient of $\mmod\ce$ by the Serre subcategory 
$\text{eff}(\ce)$, see also
Krause~\cite[Theorem~2.2]{Krause15}. In symbols Auslander's formula 
is
\[
\frac{\mmod\ce}{\text{eff}(\ce)}\eq\ce\ .
\]

Krause~\cite[Corollary~3.2]{Krause15} goes on
to give a derived category 
version
of Auslander's formula, in the derived category the formula becomes
\[
\frac{\D^b\big(\mmod\ce\big)}{\D^b_{\text{eff}(\ce)}
\big(\mmod\ce\big)}\eq
\D^b(\ce)\ .
\]
Here $\D^b_{\text{eff}(\ce)}\big(\mmod\ce\big)$ is the kernel of
the functor 
$\Lambda:\D^b\big(\mmod\ce\big)\la\D^b(\ce)$, the 
map induced on
derived categories by the 
exact functor of $\Lambda:\mmod\ce\la\ce$. Concretely
the objects
of $\D^b_{\text{eff}(\ce)}\big(\mmod\ce\big)$ are
the bounded cochain complexes in 
$\mmod\ce$ whose cohomology is in $\text{eff}(\ce)$.

Now the category $\mmod\ce$ has enough projectives, in fact
the projective objects of $\mmod\ce$ are precisely the 
essential image
of the functor $Y:\ce\la\mmod\ce$. Not only that: every
object in $\mmod\ce$ has projective dimension $\leq2$. To see this take
an object $F\in\mmod\ce$ and let 
\[\xymatrix@C+20pt{
Y(A)\ar[r]^-{Y(f)} & Y(B) \ar[r] & F\ar[r] & 0
}\]
be a finite presentation of $F$. If $K$ is the kernel in $\ce$
of the map $f:A\la B$, then the sequence
\[\xymatrix@C+20pt{
0\ar[r] & Y(K)\ar[r] &Y(A)\ar[r]^-{Y(f)} & Y(B) \ar[r] & F\ar[r] & 0
}\]
is easily seen to be exact in $\Hom\big(\ce\op,\ab\big)$, and it
exhibits a projective resolution of $F$ in the category
$\mmod\ce$ of length $\leq2$.
Thus every object in
$\D^b\big(\mmod\ce\big)$ is isomorphic to a bounded
projective resolution, and
we obtain an equivalence of triangulated categories
\[
\K^b(\ce)\quad\cong\quad\D^b\big(\mmod\ce\big)\ .
\]
The inverse image of $\D^b_{\text{eff}(\ce)}\big(\mmod\ce\big)$
under this equivalence is the category $\A^b(\ce)$ of 
Notation~\ref{N7.1}, see Krause~\cite[top of page~674]{Krause15}. 
Of course the category
$\D^b_{\text{eff}(\ce)}\big(\mmod\ce\big)$ has an obvious, standard
\tstr\ with heart $\text{eff}(\ce)$. Thus what we have really done
in Lemma~\ref{L7.3} is prove that this \tstr\ on $\A^b(\ce)$
exists for every exact category, there is no need to assume the
category $\ce$ abelian in order to produce the \tstr. 

And for an abelian category $\ce$ we have an equivalence of
categories $\A(\ce)^\heartsuit\cong\text{eff}(\ce)$.
Thus for abelian categories $\ce$, the heart of our 
new \tstr\ agrees with Auslander's old 
subcategory $\text{eff}(\ce)\subset\mmod\ce$.
\ermk

\section{The natural map 
$\D^b\big[\A(\ce)^\heartsuit\big]\la\A^b(\ce)$}
\label{S47}

Let $\ct$ be a triangulated category with a \tstr, and 
let $\ct^\heartsuit$ be the heart. Under mild hypotheses on $\ct$,
the inclusion $\ct^\heartsuit\hookrightarrow\ct^b$ can be naturally
factored as $\ct^\heartsuit\la\D^b(\ct^\heartsuit)\stackrel F\la\ct^b$,
for a triangulated functor $F:D^b(\ct^\heartsuit)\la\ct^b$.
The reader can find the most general known version of such a result
in \cite[Theorem~5.1]{Neeman91}. But in this article a more
useful version will be the older
Beilinson, Bernstein and Deligne~\cite[Proposition~3.1.10]{BeiBerDel82},
since the proof presented there comes as a quasifunctor
of model categories and hence induces a map in 
\kth.\footnote{It is natural to wonder if the functor 
$F:D^b(\ct^\heartsuit)\la\ct^b$ depends on
the choice of enhancement, and the answer is: Not much. An enhancement of
the triangulated category allows one to give it the much-weaker
structure of ``good morphisms of triangles'' as in 
\cite{Neeman91}. In order to construct the functor $F$ on all
of $\D^b(\ct^\heartsuit)$ one
needs to extend from complexes of length $n$ to complexes of length $n+1$,
and this requires constructing a mapping cone. In the presence of an 
enhancement this is choice-free, but as the reader can see
in the proof of \cite[Theorem~5.1]{Neeman91}, in the presence of
the new axioms the choice is rigid enough to permit us to choose the 
triangle uniquely up to canonical isomorphism.

For a similar argument, but with the detail spelt out much more completely
than in the sketchy~\cite[Section~5]{Neeman91},
the reader is referred to \cite[Appendix~E]{Rizzardo-VandenBergh-Neeman18}.}

It becomes natural to wonder when the functor $F$ is an equivalence.
The next Lemma is a slight variant
of~\cite[Proposition~3.1.16]{BeiBerDel82},
and gives a necessary and sufficient condition. A further variant
of~\cite[Proposition~3.1.16]{BeiBerDel82} may be found in
Chen, Han and Zhou~\cite[Theorem~2.9]{Chen-Han-Zhou19}.

\lem{L47.1}
Let $\ct$ be a triangulated category with a \tstr, let
$\ct^\heartsuit$ be the heart of the \tstr, and let 
$F:\D^b(\ct^\heartsuit)\la\ct^b$ be the natural map. The functor $F$ is an
equivalence of categories if and only 
if\,\footnote{The reader might wish to compare the necessary
and sufficient condition above with Lurie's notion of \emph{0-complicial.}
See~\cite[Definition~C.5.3.1 and Proposition~C.5.3.2]{Lurieundercon}.} every object $t\in\ct^{\leq0}\cap\ct^b$
admits a triangle $a\la t\la b$ with $a\in\ct^\heartsuit$ and 
$b\in\ct^{\leq-1}$.
\elem

\prf
Let us start with the necessity: if the functor $F$ is an equivalence then
it suffices to produce the triangle in the category $\D^b(\ct^\heartsuit)$.
The object $t\in\D^b(\ct^\heartsuit)^{\leq0}$ is isomorphic to a cochain
complex
\[\xymatrix{
\cdots\ar[r] &T^{-4}\ar[r] &T^{-3}\ar[r] &T^{-2}\ar[r] &T^{-1}\ar[r] &T^{0}
\ar[r] &0\ar[r]&\cdots
}\] 
with $T^i\in\ct^\heartsuit$, and the cochain maps
\[\xymatrix{
\cdots\ar[r] &0\ar[r]\ar[d] &0\ar[r]\ar[d] 
&0\ar[r]\ar[d] &0\ar[r]\ar[d] &T^{0}
\ar[r]\ar@{=}[d] &0\ar[r]\ar[d]&\cdots\\
\cdots\ar[r] &T^{-4}\ar[r]\ar@{=}[d] &T^{-3}\ar[r]\ar@{=}[d] 
&T^{-2}\ar[r]\ar@{=}[d] &T^{-1}\ar[r]\ar@{=}[d] &T^{0}
\ar[r]\ar[d] &0\ar[r]\ar[d]&\cdots\\
\cdots\ar[r] &T^{-4}\ar[r] &T^{-3}\ar[r] &T^{-2}\ar[r] &T^{-1}\ar[r] &0
\ar[r] &0\ar[r]&\cdots
}\]
produce the desired triangle $a\la t\la b$.

Now for the sufficiency: we assume the existence of triangles 
$a\la t\la b$ as in the
Lemma, and need to prove the functor $F$ an equivalence.

Let $n\geq1$ be an integer, let $A$ and $B$ be objects of $\ct^\heartsuit$,
and choose any morphism $f:A\la\T^nB$. Form the triangle 
$\T^{n-1}B\la t\la A\la \T^nB$. The object $t$ belongs to
$\ct^{\leq0}\cap\ct^b$, and by hypothesis there exists a triangle
$a\la t\la b$ with $a\in\ct^\heartsuit$ and $b\in\ct^{\leq-1}$.
Now let $H:\ct\la\ct^\heartsuit$ be the usual homological
functor.
The exact sequence $H^0(a)\la H^0(t)\la H^0(b)=0$
tells us that $a=H^0(a)\la H^0(t)$ is an epimorphism in
$\ct^\heartsuit$, while the exact sequence
$H^0(t)\la H^0(A)\la H^0(\T^nB)=0$ says that
$H^0(t)\la H^0(A)=A$ is also an epimorphism. We conclude
that the composite $a\la t\la A$ is an epimorphism in $\ct^\heartsuit$,
and the composite $a\la t\la A\stackrel f\la\T B$ obviously
vanishes.
Since the epimorphism $a\la A$ can be constructed for
every $f:A\la\T^nB$,  
Beilinson, Bernstein and Deligne~\cite[Proposition~3.1.16]{BeiBerDel82} 
teaches us that $F$ must be an equivalence of categories.
\eprf

We will soon be  applying Lemma~\ref{L47.1} to the case
of $F:\D^b\big[\A(\ce)^\heartsuit\big]\la\A^b(\ce)$.
Before proceeding we quickly recall

\rmd{R47.95}
Let $\ce$ be an idempotent-complete
exact category. Our notion of ``acyclic complexes''
is designed in such a way that they go to acyclic complexes
under \emph{every exact embedding} of $\ce$ as a
subcategory of an abelian category. There exists an
exact embedding
$i:\ce\la\ca$, with $\ca$ abelian,
and such that
\be
\item
The functor $i$ reflects admissible short exact sequences.
\item
A morphism $f:x\la y$ is an admissible monomorphism in $\ce$ if
and only if $i(f)$ is a monomorphism in $\ca$.
\item
A morphism $f:x\la y$ is an admissible epimorphism in $\ce$ if
and only if $i(f)$ is a epimorphism in $\ca$.
\ee
Thomason and Trobaugh~\cite[Lemma~A.7.15]{ThomTro} proves that
the Gabriel-Quillen embedding $i:\ce\la\ca$ satisfies not only (i)
but also (iii). Dually we obtain an embedding $i':\ce\la\cb$ satisfying
(i) and (ii). To achieve (i), (ii) and (iii) we take the 
embedding of $\ce$ into 
$\ca\times\cb$.

For an embedding $i:\ce\la\ca$ satisfying (i), (ii) and (iii),
if a bounded cochain complex $T\in\D^b(\ce)$ has the property 
that $i(T)$ is acyclic outside an interval $[a,b]$, then $T$ is
isomorphic in $\D^b(\ce)$ to a cochain complex
\[\xymatrix{
0\ar[r] & T^a\ar[r] & T^{a+1}\ar[r] & \cdots\ar[r] &
T^{b-1}\ar[r] & T^b\ar[r] & 0
}\]
which vanishes outside the interval $[a,b]$.
\ermd

\dfn{D47.2.5}
An idempotent-complete exact category $\ce$ is called \emph{hereditary}
if in the category $\D^b(\ce)$
the maps $E\la\T^n F$ vanish whenever $E,F\in\ce$ and $n\geq2$.
\edfn

\pro{P47.3}
Let $\ce$ be an idempotent-complete exact category, let $\A(\ce)$ be the 
homotopy category of acyclic
complexes as in Notation~\ref{N7.1}, and let the \tstr\ on
$\A(\ce)$ be as in Lemma~\ref{L7.3}.
Let $F:\D^b\big[\A(\ce)^\heartsuit\big]\la\A^b(\ce)$ be the
natural functor from the bounded derived category of the heart
to the bounded part of the \tstr\ on $\A(\ce)$, the 
subcategory $\big[\A(\ce)\big]^b=\A^b(\ce)$.

Then the functor $F$ is an equivalence if and only if the category $\ce$
is hereditary.
\epro

\rmk{R47.4}
The term ``hereditary'' goes back to Cartan and 
Eilenberg~\cite[Section~I.5]{Cartan-Eilenberg56}, where a ring $R$ is 
called hereditary if every $R$--module 
has projective dimension $\leq1$. The
rationale was that submodules of projective modules inherit the property
of being projective. In Definition~1.2
of Helmut Lenzig's 1964 PhD thesis, an abelian category is 
declared ``hereditary'' if it has global dimension $\leq1$. In 
Definition~\ref{D47.2.5} we simply extended the 
classical term to cover 
exact categories.
\ermk

\prf
Suppose the functor $F$ is an equivalence, let $n\geq2$ be an
integer and choose a morphism $f:A\la \T^nB$
in $\D^b(\ce)$. We need to show that $f$ vanishes.
Complete $f$ in $\D^b(\ce)$
to a triangle $A\stackrel f\la\T^nB\la T$. 
The object $T$ is such that every exact embedding
$i:\ce\la\ca$ takes $T$ to
a complex acyclic outside the interval $[-n,-1]$, and
Reminder~\ref{R47.95} tells us that $T$ is isomorphic
in $\D^b(\ce)$ to a cochain complex
\[\xymatrix{
\cdots\ar[r]& 0\ar[r] & T^{-n}\ar[r] &T^{-n+1}\ar[r]& \cdots\ar[r] & T^{-2} \ar[r] & 
   T^{-1} \ar[r] & 0\ar[r] &\cdots
}\]
And the triangle $\T^nB\la T\la\T A$ can be represented by cochain maps
\[\xymatrix{
\cdots\ar[r]& 0\ar[r] & B\ar[r]\ar[d] &0\ar[r]\ar[d]& \cdots\ar[r] 
& 0\ar[d] \ar[r] & 
  0 \ar[r]\ar[d] & 0\ar[r] &\cdots\\
\cdots\ar[r]& 0\ar[r] & T^{-n}\ar[d]\ar[r] &T^{-n+1}\ar[d]\ar[r]& \cdots
\ar[r] & T^{-2}\ar[d] \ar[r] & 
   T^{-1}\ar[d] \ar[r] & 0\ar[r] &\cdots\\
\cdots\ar[r]& 0\ar[r] & 0\ar[r] &0\ar[r]& \cdots\ar[r] & 0 \ar[r] & 
   A \ar[r] & 0\ar[r] &\cdots
}\]
This being a triangle means that
the sequence 
\[\xymatrix@R-20pt{
\cdots \ar[r] &0\ar[r] &
B\ar[r] & T^{-n}\ar[r] &\cdots \ar[r] & 
   T^{-1}\ar[r] &A \ar[r] & 0\ar[r]& \cdots 
}\]
is an object $t\in\A(\ce)^{\leq0}\cap\A^b(\ce)$. Since we are assuming that
$F$ is an equivalence, there exists in $\A(\ce)$ a triangle
$a\stackrel\ph\la t\la b$, with $a\in\A(\ce)^\heartsuit$ and
$b\in\A(\ce)^{\leq-1}$.

The morphism $\ph:a\la t$ in the category $\A(\ce)$
may be represented
by a cochain map
\[\xymatrix{
\cdots \ar[r] & 0\ar[r]\ar[d] & A^{-2} \ar[r]\ar[d] & 
   A^{-1} \ar[r]\ar[d] & A^{0} \ar[r]\ar[d] &0\ar[r]\ar[d]& \cdots \\
\cdots \ar[r] & T^{-3}\ar[r] & T^{-2} \ar[r] & 
   T^{-1} \ar[r] & A \ar[r] &0\ar[r]& \cdots 
}\]
The distinguished triangle $a\stackrel\ph\la t\la b$ tells us that $b$ is 
homotopy equivalent to $\text{Cone}(\ph)$, the mapping cone on the cochain 
map $\ph$. Requiring that $b\cong\text{Cone}(\ph)$ should belong to
$\A(\ce)^{\leq-1}$ amounts to saying that the morphism 
$A^0\oplus T^{-1}\la A$ must be a split epimorphism in $\ce$.

Now the morphism $a\la t$ is isomorphic in $\A(\ce)$ to the cochain map
\[\xymatrix{
\cdots \ar[r] & 0\ar[r]\ar[d] & A^{-2} \ar[r]\ar[d] & 
   A^{-1}\oplus T^{-1} \ar[r]\ar[d] & A^{0}\oplus T^{-1} \ar[r]\ar[d] &0\ar[r]\ar[d]& \cdots \\
\cdots \ar[r] & T^{-3}\ar[r] & T^{-2} \ar[r] & 
   T^{-1} \ar[r] & A \ar[r] &0\ar[r]& \cdots 
}\]
and thus, without changing the isomorphism class of
the map $\ph:a\la t$ in the category $\A(\ce)$, we may assume that the
cochain map is such that $A^0\la A$ is a split epimorphism. Choose
a splitting, meaning choose a morphism $g:A\la A^0$ such that the composite
$A\la A^0\la A$ is the identity. Now form in $\ce$ the pullback square
\[\xymatrix{
B^{-1} \ar[r]\ar[d] & A\ar[d]^g \\
A^{-1} \ar[r] & A^{0}
}\]
Then the composite
\[\xymatrix{
\cdots \ar[r] & 0\ar[r]\ar[d] & A^{-2} \ar[r]\ar@{=}[d] & 
   B^{-1} \ar[r]\ar[d] & A \ar[r]\ar[d]^g &0\ar[r]\ar[d]& \cdots \\
\cdots \ar[r] & 0\ar[r]\ar[d] & A^{-2} \ar[r]\ar[d] & 
   A^{-1} \ar[r]\ar[d] & A^{0} \ar[r]\ar[d] &0\ar[r]\ar[d]& \cdots \\
\cdots \ar[r] & T^{-3}\ar[r] & T^{-2} \ar[r] & 
   T^{-1} \ar[r] & A \ar[r] &0\ar[r]& \cdots 
}\]
is a cochain map between acyclic complexes
\[\xymatrix{
\cdots \ar[r] & 0\ar[r]\ar[d] & A^{-2} \ar[r]\ar[d] & 
   B^{-1} \ar[r]\ar[d] & A \ar[r]\ar@{=}[d] &0\ar[r]\ar[d]& \cdots \\
\cdots \ar[r] & T^{-3}\ar[r] & T^{-2} \ar[r] & 
   T^{-1} \ar[r] & A \ar[r] &0\ar[r]& \cdots 
}\]
Coming back to the category $\D^b(\ce)$: the cochain maps
\[\xymatrix{
\cdots \ar[r] &0\ar[r] & 0\ar[d]\ar[r]& \cdots\ar[r] &
 0\ar[r]\ar[d] & A^{-2} \ar[r]\ar[d] & 
   B^{-1} \ar[r]\ar[d] &0\ar[r]\ar[d]& \cdots \\
\cdots \ar[r] &0\ar[r] & T^{-n}\ar[r]\ar[d]& \cdots\ar[r] & T^{-3}\ar[d] \ar[r]
& T^{-2} \ar[d]\ar[r] & 
   T^{-1}\ar[d] \ar[r] & 0\ar[r]& \cdots \\
\cdots \ar[r] &0\ar[r] & 0\ar[r]& \cdots\ar[r] & 0 \ar[r]
&0 \ar[r] & 
   A \ar[r] & 0\ar[r]& \cdots 
}\]
can be viewed as morphisms $A\la \T^{-1}T\la A$ composing to
the identity. The triangle $\T^{-1}T\la A\stackrel f\la\T^n B\la T$
is split, and the map $f:A\la\T^nB$ that we started 
out with must vanish
in $\D^b(\ce)$.

Suppose now that in the category $\D^b(\ce)$ any morphism $A\la\T^2B$
vanishes, whenever $A,B\in\ce$. Any object in $t\in\A(\ce)^{\leq0}\cap\A^b(\ce)$ 
of the
form
\[\xymatrix{
0\ar[r]  & T^{-3}\ar[r] & T^{-2} \ar[r] &
   T^{-1} \ar[r] & T^{0} \ar[r] &0 
}\]
gives rise to a morphism $T^0\la\T^2 T^{-3}$ in the category 
$\D^b(\ce)$, and it is classical that this morphism will vanish
in $\D^b(\ce)$ if and only if there is a cochain map of acyclic
complexes
\[\xymatrix{
0\ar[r]  & 0\ar[d]\ar[r] & C^{-2} \ar[r]\ar[d] &
   C^{-1} \ar[r]\ar[d] & T^{0} \ar[r]\ar@{=}[d] &0 \\
0\ar[r]  & T^{-3}\ar[r] & T^{-2} \ar[r] &
   T^{-1} \ar[r] & T^{0} \ar[r] &0 
}\]
This produces for us a morphism $a\la t$ with $a\in\A(\ce)^\heartsuit$
whose mapping cone lies in $\A(\ce)^{\leq-1}$.

If $t\in\A(\ce)^{\leq0}\cap\A^b(\ce)$ is more general, meaning of the form
\[\xymatrix{
0\ar[r]  & T^{-n}\ar[r] & T^{-n+1} \ar[r] &\cdots\ar[r] &
   T^{-1} \ar[r] & T^{0} \ar[r] &0 
}\]
we may apply the above to the complex
\[\xymatrix{
0\ar[r]  & K\ar[r] & T^{-2} \ar[r] &
   T^{-1} \ar[r] & T^{0} \ar[r] &0 
}\]
with $K$ the image of $T^{-3}\la T^{-2}$; since $t$ is acyclic this 
image is in $\ce$. This produces for us a cochain map
\[\xymatrix{
0\ar[r]  & 0\ar[d]\ar[r] & C^{-2} \ar[r]\ar[d] &
   C^{-1} \ar[r]\ar[d] & T^{0} \ar[r]\ar@{=}[d] &0 \\
0\ar[r]  & K\ar[r] & T^{-2} \ar[r] &
   T^{-1} \ar[r] & T^{0} \ar[r] &0 
}\]
and hence also a cochain map
\[\xymatrix{
\cdots\ar[r]  & 0\ar[d]\ar[r] & C^{-2} \ar[r]\ar[d] &
   C^{-1} \ar[r]\ar[d] & T^{0} \ar[r]\ar@{=}[d] &0 \\
\cdots\ar[r]  & T^{-3}\ar[r] & T^{-2} \ar[r] &
   T^{-1} \ar[r] & T^{0} \ar[r] &0 
}\]
which may be viewed as a morphism $\ph:a\la t$ in $\A(\ce)$, with
$a\in\A(\ce)^\heartsuit$ and such that the mapping cone of 
$\ph$ belongs to $\A(\ce)^{\leq-1}$.
\eprf

\section{The $K$--theoretic consequences}
\label{S8}

\rmd{R8.-1}
Let $\ce$ be an idempotent-complete exact categry.
In Notation~\ref{N7.1} we recalled the categories $\A^?(\ce)\subset\K^?(\ce)$
and the Verdier quotient $\D^?(\ce)=\K^?(\ce)/\A^?(\ce)$, where
$?$ is any of $b$, $+$, $-$ or the empty restriction.
Now a special case is the exact category $\ce^\oplus$. This means
we take any idempotent-complete additive category $\ce$, and give
it the exact structure where the admissible exact sequences are 
the split short exact sequences.

Specializing to the case
of $\ce^\oplus$ the 
general definitions of Notation~\ref{N7.1},
an acyclic complex $E^*\in\A(\ce^\oplus)$ is a cochain complex
\[\xymatrix@C+20pt{
\cdots \ar[r]^-{\partial^{i-2}} & E^{i-1}\ar[r]^-{\partial^{i-1}} &
E^{i}\ar[r]^-{\partial^{i}} &
E^{i+1}\ar[r]^-{\partial^{i+1}} &\cdots
}\] 
where there exist split short exact sequences
\[\xymatrix@C+20pt{
0 \ar[r] & K^{i}\ar[r]^-{\alpha^{i}} &
E^{i}\ar[r]^-{\beta^{i}} &
K^{i+1}\ar[r] &0
}\] 
such that $\partial^i=\alpha^{i+1}\circ\beta^i$. 
This makes $E^i\cong K^i\oplus K^{i+1}$, and the 
complex $E^*$ can be decomposed as a direct sum of complexes
\[\xymatrix@C+20pt{
\cdot\ar[r] &0 \ar[r] & K^{i}\ar[r]^-{\id} &
K^{i}\ar[r] &0\ar[r] &\cdots
}\]
which vanish except in degrees $(i-1)$ and $i$. Hence all
objects in $\A^?(\ce^\oplus)$ are contractible, and are
isomorphic to zero in the homotopy category
$\K^?(\ce^\oplus)=\K^?(\ce)$. This makes the Verdier quotient
\[\D^?(\ce^\oplus)\eq\K^?(\ce^\oplus)/\A^?(\ce^\oplus)\eq\K^?(\ce)/0\eq\K^?(\ce).\] 
\ermd

\rmd{R8.-3}
Now it's time to move on to the {\it K--}theoretic consequences,
which means we need to consider model categories as well as
the associated triangulated categories.
In the remainder of this section we follow the conventions of 
Schlichting~\cite{Schlichting06}. Thus $\cm$ will be a category 
of models, $\D:\cm\la\ct$ will be a functor from $\cm$ to the category
$\ct$ of small triangulated categories, and $\KK$ will be a functor
from $\cm$ to spectra. And we will assume that if $M'\la M\la M''$
is an exact sequence in $\cm$
then
\[\xymatrix@C+20pt{
\KK(M')\ar[r] &\KK(M)\ar[r] &\KK(M'')
}\]
is a homotopy fibration. Recall: 
the sequence  $M'\la M\la M''$ is declared to be exact in $\cm$ if
the categories $\D(\cm')$, $\D(\cm)$ and $\D(\cm'')$ are
all idempotent-complete, and
(1) the functor $\D(M')\la\D(M)$ is fully 
faithful, (2) the composite $\D(M')\la\D(M)\la\D(M'')$ vanishes, and
(3) the natural map $\D(M)/\D(M')\la\D(M'')$ is fully faithful, with
$\D(\cm'')$ being the idempotent-completion of the essential image
of $\D(M)/\D(M')$.
\ermd

The next result was presaged in 
Schlichting~\cite[Proposition~2 in Section~11]{Schlichting06}. 
What's incorrect about
Schlicting's proof of his old 
proposition could be rephrased as saying that
the natural functor
\[\xymatrix@C+30pt{
\D^b\big[\A(\ce)^\heartsuit\big] \ar[r] &
\A^b(\ce)
}\]
need not in general be an equivalence; see
Proposition~\ref{P47.3}.

We begin with the easy

\pro{P8.1}
Let $\ce$ be an idempotent-complete exact category, and let
$\ce^\oplus$ be the category $\ce$ but with the split exact structure.
Then there is a homotopy fibration of non-connective
\kth-spectra
\[\xymatrix@C+30pt{
\KK(M')\ar[r] & \KK\big(\ce^\oplus\big) \ar[r] &\KK(\ce)
}\]
where $M'\in\cm$ satisfies $\D(M')=\A^b(\ce)$.
\epro

\prf
Let $\cm$ be the category of biWaldhausen complicial categories
as in Thomason and Trobaugh~\cite{ThomTro}.
We consider the sequence $M'\la M\la M''$ in $\cm$ where,
in the notation of Schlichting~\cite{Schlichting06}---which is
in conflict with our notation---one would write
$M'=\A^b(\ce)$, $M=\mathbf{Ch}^b(\ce)$ and
$M''=\big(\mathbf{Ch}^b(\ce,\A^b(\ce))\big)$. The conflict of
notation is that
in \cite{Schlichting06} $\A^b(\ce)$ is an object of $\cm$,
whereas our $\A^b(\ce)$ is what in Schlichting's notation would be
$\D\big[\A^b(\ce)\big]$. Our notation, which [as we have said] clashes
with Schlichting, follows Krause~\cite{Krause15}. It is impossible to
choose a notation which agrees with every previous article in the 
literature.

 We remind the reader:
Schlichting's notation means that the 
objects and morphisms of $M,M''$ are the same, both categories have
for objects the bounded 
complexes of objects in $\ce$ and the morphisms 
are the cochain maps. The weak equivalences in $M$ are
the homotopy equivalences, whereas the weak equivalences in $M''$
are the cochain maps whose mapping cones are acyclic. The objects
of $M'$ are the acyclic complexes, and the weak equivalences are
as in $M$. The natural functor $\D:\cm\la\ct$ takes the 
sequence $M'\la M\la M''$  
to the sequence of triangulated
categories $\A^b(\ce)\la\K^b(\ce)\la\D^b(\ce)$, where
$\A^b(\ce)$ is to be understood in our notation, it
is a triangulated category. 

The categories $\A^b(\ce)$,
$\K^b(\ce)$ and $\D^b(\ce)$ are all known to be idempotent-complete:
in the case of $\D^b(\ce)$ this is
by \cite[Theorem~2.8]{Balmer-Schlichting}. For
$\K^b(\ce)$, Reminder~\ref{R8.-1} tells us that
that $\K^b(\ce)=\D^b(\ce^\oplus)$, reducing us to the previous case.
And for $\A^b(\ce)\subset\K^b(\ce)$ this is because
we know $\A^b(\ce)$ to be a thick subcategory of
the idempotent-complete triangulated category
$\K^b(\ce)$, we already mentioned the thickness of
$\A^b(\ce)$ as a subcategory of $\K^b(\ce)$ in
Notation~\ref{N7.1}.
Hence the sequence $M'\la M\la M''$ is exact in $\cm$. 
Therefore the functor
$\KK$ takes  $M'\la M\la M''$ to 
a homotopy fibration.
In this homotopy fibration we have that $\KK(M'')=\KK(\ce)$,
just because $\D(M'')=\D^b(\ce)$. And $\KK(M)=\KK(\ce^\oplus)$,
on the grounds that $\D(M)=\K^b(\ce)=\D^b(\ce^\oplus)$.
Thus our homotopy fibration becomes 
$\KK(M')\la\KK\big(\ce^\oplus\big)\la\KK(\ce)$.
\eprf

\rmk{R8.3}
Proposition~\ref{P8.1} was straightforward, but in combination with
Lemma~\ref{L7.3} it becomes remarkable. The homotopy
fiber on the map $\KK\big(\ce^\oplus\big)\la\KK(\ce)$ is 
identified with $\KK(M')$, and $\D(M')\cong\A^b(\ce)$ has 
a bounded \tstr. 

Antieau, Gepner and 
Heller~\cite[Conjecture~B on page~244]{Antieau-Gepner-Heller19},
if true, would imply 
that $\KK_{-n}(M')=0$ for all $n>0$, and we would deduce
that the natural map $\KK_{-n}\big(\ce^\oplus\big)\la\KK_{-n}(\ce)$
would have to be an isomorphism for all $n>0$. But this will be shown
to be false, see Example~\ref{E8.5} below.
As it happens in Example~\ref{E8.5} the 
exact category $\ce$ will be hereditary, and Proposition~\ref{P47.3}
informs us that the natural map
$F:\D^b\big[\A(\ce)^\heartsuit\big]\la\A^b(\ce)$ is an
equivalence of categories. From
Beilinson, Bernstein and 
Deligne~\cite[proof of Proposition~3.1.10]{BeiBerDel82}
we know that the map $F$ can be realized as
$\D(f)$ for some suitable morphism in $\cm$, and hence 
induces an isomorphism in \kth. It follows
that $\KK_{-n}\big[\A(\ce)^\heartsuit\big]$ is also
nonzero for some $n>0$, contradicting 
\cite[Conjecture~A]{Antieau-Gepner-Heller19},
which is a restatement of an older conjecture due to 
Schlichting~\cite[Conjecture~1 of Section~10]{Schlichting06}.

There exists an abelian category $\A(\ce)^\heartsuit$ 
with non-vanishing negative \kth.
\ermk

\section{Vector bundles on nodal curves}
\label{S470}

We begin by recalling classical facts about smooth projective curves.

\rmd{R47.1}
Fix an algebraically closed field $k$, and 
let $X$ be a smooth, projective curve over $k$.
We allow $X$ to have more than one connected component.
A vector bundle on $X$ will mean a locally free sheaf,
locally of finite rank. Of course: the rank may depend
on the connected component we're at.

Each vector bundle $\cv$ on $X$ gives rise to three
continuous functions 
\[
\text{\rm rank}:X\la\nn\,,\qquad\text{\rm degree}:X\la\zz\,,\qquad
\text{\rm slope}:X\la\qq\,.
\]
These continuous functions assign a number to each connected
component of $X$; for the rank this number is a positive integer,
the degree may be any integer, and the slope is a rational number.
The rank is obvious. The degree of a line bundle,
on a connected component $X_i\subset X$,  is the usual degree---the 
number of zeros minus the number of poles of a rational
section. For a vector bundle
$\cv$ of rank $n$ the degree of $\cv$ is defined to be the degree 
of the line bundle $\wedge^n\cv$, and of course it will depend
on the component. And the slope---whose early history will be 
recalled in Remark~\ref{R47.2}---is defined by the formula
\[
\text{\rm slope}\left(\cv\right)\eq\frac{\text{\rm degree}\left(\cv\right)}
{\text{\rm rank}\left(\cv\right)}\ ,
\]
and we repeat that this fomula yields a rational number for each
component of $X$.

We remind the reader that, if $f:\cv\la\cv'$ is an injective map
of vector bundles of equal rank, then 
$\text{\rm degree}(\cv)\leq \text{\rm degree}(\cv')$ with
equality if and only if $f$ is an isomorphism. 
Since for us this fact is key
we quickly recall
the argument:
assume $X$ connected, and let $n=\text{\rm rank}(\cv)=\text{\rm rank}(\cv')$.
The determinant $\wedge^nf:\wedge^n\cv\la\wedge^n\cv'$ 
is an injective map between line bundles, and may be viewed
as a regular, nonzero section of the line bundle
$\cl=(\wedge^n\cv')\otimes(\wedge^n\cv)^{-1}$.
Being a regular section of $\cl$, the determinant
$\wedge^nf$ has no poles. The degree of $\cl$ can be computed as
the number of zeros of $\wedge^nf$ and must be non-negative---therefore
$\text{\rm degree}(\cv')-\text{\rm degree}(\cv)=\text{\rm degree}(\cl)\geq0$.
Equality is equivalent to $\wedge^nf$ having no zeros, which happens if and
only if $\wedge^nf:\wedge^n\cv\la\wedge^n\cv'$ is an isomorphism. But this 
is equivalent to $f:\cv\la\cv'$ being an isomorphism.

And for us the most important consequence is: 
any injective endomorphism of a vector bundle $\cv$ must be an 
isomorphism.
\ermd

\rmk{R47.2}
Since I've been asked about the history:
it's traditional
to assume $X$ irreducible, which we will do in this
Remark. Identifying both
$H^0(X)$ and $H^2(X)$ with
$\zz$, we have that
$\text{\rm rank}(\cv)=\text{ch}_0^{}(\cv)$ and
$\text{\rm degree}(\cv)=\text{ch}_1^{}(\cv)$ are the zeroth and
first Chern characters of $\cv$ and go back a long way. The quotient 
$\mu(\cv)=\text{\rm degree}(\cv)/\text{\rm rank}(\cv)$
was first explicitly introduced in print
in 1969 by Narasimhan and 
Ramanan~\cite[beginning of Section~2]{Narasimhan-Ramanan69}.
The name ``slope'' for the rational number $\mu(\cv)$
came later---Ramanan tells me that it was coined by Quillen,
but the first occurrence I can find in print is 
in the 1977 article by
Shatz~\cite[Section~2]{Shatz77}.
Of course in some sense it all goes back 
Mumford's 1962 ICM talk~\cite[Definition on page~529]{Mumford63},
which discusses 
the formation of the moduli space of stable vector
bundles. 
\ermk

\dis{D47.3}
Now it's time to move on to singular curves---but for 
simplicity the only singularities
we allow are simple nodes.\footnote{The argument we
are about to sketch hinges on the relation between
a vector bundle $\cv$, on a singular curve $Y$, and
its pullback $\pi^*\cv$ via the normalization map
$\pi:X\la Y$. By confining attention to simple
double points we make our life easier---the map 
$\pi:X\la Y$ is especially straightforward then. But the argument
extends, and the reader wishing to treat
the general case is referred to
Serre~\cite[Chapter~IV, Section~1]{Serre59}
for a discussion of the normalization map
$\pi:X\la Y$ for arbitrary singular curves $Y$.}
Let $Y$ be a nodal projective curve, and let $\pi:X\la Y$ be
the normalization. Then $X$ is a smooth projective curve as
in Reminder~\ref{R47.1}. A vector bundle $\cv$ on $Y$, which is a locally
free sheaf on $Y$ locally of finite rank, pulls back to a vector
bundle $\pi^*\cv$ on $X$. And by Reminder~\ref{R47.1}
the vector bundle $\pi^*\cv$ has associated to it
rank, degree and slope functions.

Let $\{p_1^{},p_2^{}\ldots,p_n^{}\}$ be the singular points of $Y$. Then
each $p_i^{}$ has two 
distinct inverse images in $X$; let us call them $p'_i$ and 
$p''_i$. Thus for each $i$ we have a commutative diagram of schemes
\[\xymatrix@C+30pt@R+10pt{
\spec k\ar[r]^-{\alpha_i}\ar[d]_-{\beta_i}\ar[dr]^{\gamma_i} & X\ar[d]^\pi \\
X\ar[r]_-\pi & Y
}\]
where the image of $\alpha_i$ is $p'_i$ and the image of $\beta_i$ is 
$p''_i$. The composite $\gamma_i=\pi\alpha_i=\pi\beta_i$
has image $p_i=\pi(p'_i)=\pi(p''_i)$. 
For any coherent sheaf $\cw$ on $X$ we have units of adjunction
$\eta(\alpha_i):\cw\la \alpha_{i*}\alpha_i^*\cw$ and
$\eta(\beta_i):\cw\la \beta_{i*}\beta_i^*\cw$. Now let $\cv$ be a
coherent sheaf on $Y$. There is the unit of adjunction
$\eta(\pi):\cv\la\pi_*\pi^*\cv$, and we may combine with 
the units of adjunction above to obtain a commutative square
\[\xymatrix@C+70pt{
\gamma_{i*}\gamma_i^*\cv=\pi_*\alpha_{i*}\alpha_i^*\pi^*\cv=\pi_*\beta_{i*}\beta_i^*\pi^*\cv&\pi_*\pi^*\cv\ar[l]_-{\eta(\alpha_i)}  \\
\pi_*\pi^*\cv\ar[u]^{\eta(\beta_i)} & \ar[l]^-{\eta(\pi)}\ar[u]_{\eta(\pi)}\cv
}\]
And it is a classical fact that, for 
$\cv$ a vector bundle on $Y$, the
sequence
\[\xymatrix@C+15pt{
0\ar[r]& \cv\ar[r]^-{\eta(\pi)} &\pi_*\pi^*\cv
\ar[rrr]^-{\oplus_{i=1}^n\big[\eta(\alpha_i)-\eta(\beta_i)\big]} 
&&& \ds\bigoplus_{i=1}^n
\gamma_{i*}\gamma_i^*\cv
\ar[r] & 0
}\]
is an exact sequence of sheaves on $Y$. Moreover: this sequence
can be used to construct vector bundles on $Y$. A vector bundle $\cv$ on
$Y$ is uniquely determined by the vector bundle $\cw=\pi^*\cv$ on
$X$, together with the isomorphisms $\alpha_i^*\cw\cong\beta_i^*\cw$
that arise from the canonical isomorphism 
$\alpha_i^*\pi^*\cv\cong\beta_i^*\pi^*\cv$. Concretely:
given any vector bundle $\cw$ on $X$ and, for each $i$, an isomorphism
$\alpha_i^*\cw\cong\beta_i^*\cw$,
we define $\cv$ to be the kernel of the map of
sheaves
\[\xymatrix@C+10pt{
 \pi_*\cw
\ar[rrr]^-{\oplus_{i=1}^n\big[\eta(\alpha_i)-\eta(\beta_i)\big]} 
&&& \ds\bigoplus_{i=1}^n
\big[\gamma_{i_*}\alpha_i^*\cw\cong \gamma_{i*}\beta_i^*\cw\big]
}\]
where the isomorphism $\gamma_{i*}\alpha_i^*\cw\cong \gamma_{i*}\beta_i^*\cw$
is by applying the functor $\gamma_{i*}$
to the chosen isomorphism $\alpha_i^*\cw\cong\beta_i^*\cw$.
By construction $\cv$ is a coherent sheaf on $Y$. 
And checking that $\cv$ is 
a vector bundle on $Y$, with the natural map $\pi^*\cv\la\cw$
an isomorphism, 
is local in $Y$ in the flat topology. Hence we may do it
separately
at each singular point $p_i\in Y$, 
and simplify the argument by first completing at $p_i$.
This we leave to the reader.

Moreover: morphisms of vector bundles $f:\cv\la\cv'$ on $Y$ are uniquely
determined by the ``descent data'' above. The morphism $f$ gives
rise to a map of short exact sequences
\[\xymatrix@R+40pt{
0\ar[r]& \cv\ar[r]^-{\eta(\pi)} \ar[d]_f&\pi_*\pi^*\cv\ar[d]^{\pi_*\pi^*f}
\ar[rrr]^-{\oplus_{i=1}^n\big[\eta(\alpha_i)-\eta(\beta_i)\big]} 
&&& \ds\bigoplus_{i=1}^n
\big[\gamma_{i*}\alpha_i^*\pi^*\cv=\gamma_{i*}\beta_i^*\pi^*\cv\big]
\ar[r]\ar[d]|{\oplus_{i=1}^n
\big[\gamma_{i*}\alpha_i^*\pi^*f=\gamma_{i*}\beta_i^*\pi^*f\big]}  & 0 \\
0\ar[r]& \cv'\ar[r]^-{\eta(\pi)} &\pi_*\pi^*\cv'
\ar[rrr]^-{\oplus_{i=1}^n\big[\eta(\alpha_i)-\eta(\beta_i)\big]} 
&&& \ds\bigoplus_{i=1}^n
\big[\gamma_{i*}\alpha_i^*\pi^*\cv'=\gamma_{i*}\beta_i^*\pi^*\cv'\big]
\ar[r]& 0
}\]
and hence the commutative square on the right uniquely determines the 
vertical arrow on the left.
\edis

\rmk{R47.314}
Just to clarify: everything in Discussion~\ref{D47.3} is known. See
Theorem~1.3 on page 18 of Igor Burban's 2003 PhD thesis for a complete
and thorough treatment, and Burban and 
Kreussler~\cite[Theorem~5.1.4]{Burban-Kreussler12}
for a published account. For a partial list of articles
on related topics, going back to the 1960s, see
Burban and Drozd~\cite[beginning of Chapter~2]{Burban-Drozd17}.
\ermk

And now the time has come to prove something.

\lem{L47.5}
Let $k$ be an algebraically closed field, let $Y$ be a projective curve over
$k$ with only simple nodes, and  
let $\cv$ be a vector bundle on $Y$. Then any self-map $f:\cv\la\cv$
leads to a canonical decomposition $\cv=\cv'\oplus\cv''$
such that 
\be
\item
$f$ decomposes as $f=f'\oplus f''$ for maps $f':\cv'\la\cv'$
and $f'':\cv''\la\cv''$.
\item
The map $f'$ is an isomorphism while the map $f''$ is nilpotent.
\setcounter{enumiv}{\value{enumi}}
\ee
\elem

\rmk{R47.6}
The beginning of the proof below, which will treat the case where $Y$ is smooth, has
similarities with the proof of Fitting's Lemma, see 
Jacobson~\cite[pp.~113-114]{Jacobson89}. Since the category of 
coherent sheaves on $Y$ isn't Artinian the arguments aren't
identical; nevertheless the reader might wish
to compare the two proofs.
\ermk

\prf
We first treat the case where $Y$ is smooth. The coherent sheaves
$\Ker(f^n)$ are an increasing sequence of subsheaves of the coherent sheaf
$\cv$, hence must stabilize. There exists an integer $N\gg0$ such that,
for all $n\geq N$, the inclusion $\Ker(f^n)\subset\Ker(f^{n+1})$ is
an isomorphism.

Now: for each integer $n>0$ we have a short exact sequence of coherent
sheaves
\[\xymatrix{
0\ar[r] & \Ker(f^n) \ar[r]^-{\psi_n} & \Ker(f^{2n})\ar[r] & \Ker(f^n)\cap\Ima(f^n)\ar[r] & 0
}\]
and if $n\geq N$ the map $\psi_n$ is an isomorphism. We
conclude that, as long as $n\geq N$, we must have
$\Ker(f^n)\cap\Ima(f^n)=0$.  

Next note that $\Ker(f^n)$ and $\Ima(f^n)$ are coherent subsheaves of
$\cv$. They must be torsion-free, and torsion-free coherent sheaves
on a smooth curve are vector bundles. Choose
any $n\geq N$ and 
consider the natural composable morphisms
\[\xymatrix{
 \Ima(f^n) \ar[r] & \cv\ar[r] & \Ima(f^n)
}\]
The kernel of the composite is zero, 
since the map $\cv\la\Ima(f^n)$ has kernel
$\Ker(f^n)$ which intersects $\Ima(f^n)\subset\cv$ trivially.
Thus the composite is an injective endomorphism of the vector bundle 
$\Ima(f^n)$ and
must be an isomorphism. We conclude that $\Ima(f^n)$ is a direct
summand of $\cv$, more precisely we learn that the inclusion
$\Ima(f^n)\la\cv$ provides a splitting of the
short exact sequence
\[\xymatrix{
0\ar[r] & \Ker(f^n) \ar[r] & \cv\ar[r] & \Ima(f^n)\ar[r] & 0\ .
}\]
This gives us our canonical decomposition 
$\cv=\Ima(f^n)\oplus\Ker(f^n)$. Obviously the map
$f:\cv\la\cv$ takes $\Ima(f^n)\subset\cv$ to itself and takes
$\Ker(f^n)\subset\cv$ to itself, that is $f:\cv\la\cv$ 
decomposes as $f'\oplus f''$ for a unique
$f':\Ima(f^n)\la\Ima(f^n)$ and a unique $f'':\Ker(f^n)\la\Ker(f^n)$.
Obviously $f''$ is nilpotent, more precisely $(f'')^n=0$. 
And the kernel of the map $f'$ is
\[
\Ker(f)\cap\Ima(f^n)\sub\Ker(f^n)\cap\Ima(f^n)\eq0\ ,
\]
which shows that the map $f':\Ima(f^n)\la\Ima(f^n)$ is an injective
endomorphism of the vector bundle $\Ima(f^n)$ and hence an 
isomorphism.

Now for the case where $Y$ is allowed to have simple nodes as singularities.
Let $f:\cv\la\cv$ be an endomorphism of the vector bundle $\cv$,
and let $\pi:X\la Y$ be the normalization of $Y$. Then 
$\pi^*f:\pi^*\cv\la\pi^*\cv$ is an endomorphism of the vector bundle 
$\pi^*\cv$ on $X$, and by the above there exists a decomposition
of $\pi^*\cv$ as $\pi^*\cv=\cw'\oplus\cw''$ such that
\be
\setcounter{enumi}{\value{enumiv}}
\item
There exists an integer $n\gg0$ with $\cw'=\Ima(\pi^*f^n)$
and $\cw''=\Ker(\pi^*f^n)$.
\item
The map 
$\pi^* f:\cw'\oplus\cw''\la\cw'\oplus\cw''$ is equal to 
$\ph'\oplus \ph''$ for 
some morphisms $\ph':\cw'\la\cw'$ and $\ph'':\cw''\la\cw''$.
\item
The map $\ph'$ is an isomorphism while the map $\ph''$
satisfies $(\ph'')^n=0$.
\setcounter{enumiv}{\value{enumi}}
\ee
With $\alpha_i:\spec k\la X$ and $\beta_i:\spec k\la X$ as
in Discussion~\ref{D47.3}, we have induced decompositions
$\alpha_i^*\pi^*\cv=\alpha_i^*\cw'\oplus\alpha_i^*\cw''$
and 
$\beta_i^*\pi^*\cv=\beta_i^*\cw'\oplus\beta_i^*\cw''$ such that
\be
\setcounter{enumi}{\value{enumiv}}
\item
$\alpha_i^*\pi^* f=\alpha_i^*\ph'\oplus\alpha_i^*\ph''$ 
and $\beta_i^*\pi^* f=\beta_i^*\ph'\oplus\beta_i^*\ph''$.
\setcounter{enumiv}{\value{enumi}}
\ee
Of course we have 
canonical isomorphisms $\rho_i:\alpha_i^*\pi^*\cv\la\beta_i^*\pi^*\cv$,
and these ismorphisms $\rho_i$ must be such that the squares below commute
\[\xymatrix@C+40pt{
\alpha_i^*\pi^* \cv \ar[r]^-{\alpha_i^*\pi^* f} 
\ar[d]_{\rho_i}& \alpha_i^*\pi^* \cv\ar[d]^{\rho_i} \\
\beta_i^*\pi^* \cv \ar[r]^-{\beta_i^*\pi^* f} & \beta_i^*\pi^* \cv
}\]
Raising the horizontal maps to the $n\mth$ power, the squares
\[\xymatrix@C+40pt{
\alpha_i^*\pi^* \cv \ar[r]^-{\alpha_i^*\pi^* f^n} 
\ar[d]_{\rho_i}& \alpha_i^*\pi^* \cv\ar[d]^{\rho_i} \\
\beta_i^*\pi^* \cv \ar[r]^-{\beta_i^* \pi^*f^n} & \beta_i^*\pi^* \cv
}\]
must also commute.
But these squares can be rewritten as
\[\xymatrix@C+70pt{
\alpha_i^*\cw'\oplus\alpha_i^*\cw'' \ar[r]^-{\big[\alpha_i^*(\ph')^n\big]\oplus0} 
\ar[d]_{\rho_i}& \alpha_i^*\cw'\oplus\alpha_i^*\cw''\ar[d]^{\rho_i} \\
\beta_i^*\cw'\oplus\beta_i^*\cw'' \ar[r]^-{\big[\beta_i^*(\ph')^n\big]\oplus0} & 
\beta_i^*\cw'\oplus\beta_i^*\cw''
}\]
The fact that $\ph'$ is an isomorphism means that so are
$\alpha_i^*(\ph')^n$ and $\beta_i^*(\ph')^n$. And the commutativity of
the square forces the map $\rho_i$ to take the kernel of the top horizontal
map to the kernel of the bottom horizontal map, and the image of the top
horizontal map to the image of the bottom horizontal map. That is: the
isomorphism 
$\rho_i:\alpha_i^*\cw'\oplus\alpha_i^*\cw''\la\beta_i^*\cw'\oplus\beta_i^*\cw''$ 
must split as the direct sum $\rho'_i\oplus\rho''_i$, for
isomorphisms 
$\rho'_i:\alpha_i^*\cw'\la\beta_i^*\cw'$ and 
$\rho''_i:\alpha_i^*\cw''\la\beta_i^*\cw''$.

And this allows us to descend to $Y$; the maps $\ph':\cw'\la\cw'$ and 
$\ph'':\cw''\la\cw''$, together with the descent data given by the
isomorphisms $\rho'_i$ and $\rho''_i$, allow us to uniquely
define vector bundles $\cv',\cv''$ on $Y$, as well as endomorphisms 
$f':\cv'\la\cv'$ and $f'':\cv''\la\cv''$, so that
$\pi^*f':\pi^*\cv'\la\pi^*\cv'$ agrees with $\ph':\cw'\la\cw'$ and
$\pi^*f'':\pi^*\cv''\la\pi^*\cv''$ agrees with $\ph'':\cw''\la\cw''$.
The uniqueness forces $f'\oplus f''$ to agree with $f:\cv\la\cv$.
And the fact that $f'$ is an isomorphism and $(f'')^n=0$ can be checked 
after pulling back to $X$.
\eprf

\rmk{R47.7}
We should perhaps spell out what we meant in Lemma~\ref{L47.5}, when
we said that the decomposition of $\cv$ as $\cv=\cv'\oplus\cv''$ is
``canonical''.

The fact that $f:\cv\la\cv$ decomposes as $f=f'\oplus f''$, with
$f':\cv'\la\cv'$ an isomorphism and $f'':\cv''\la\cv''$ nilpotent, makes
the decomposition unique. Choose an $n\gg0$ so that $(f'')^n=0$; then
$f^n:\cv\la\cv$ has kernel $\cv''$ and image $\cv'$. Thus we can write
the formulas
\[
\cv'\eq\bigcap_{n=1}^\infty\Ima(f^n)\qquad\text{ and }\qquad
\cv''\eq\bigcup_{n=1}^\infty\Ker(f^n)\ ,
\]
which show that these subbundles are canonically unique---what isn't
immediate from the formulas is that these are vector bundles and that
their direct sum is $\cv$.

Suppose we are given a commutative square of maps
of vector bundles on $Y$
\[\xymatrix{
\cv\ar[r]^f\ar[d]_\rho & \cv\ar[d]^\rho \\
\cw\ar[r]^g &\cw
}\]
Then obviously
\[
\rho\left(\bigcap_{n=1}^\infty\Ima(f^n)\right) \subset 
           \bigcap_{n=1}^\infty\Ima(g^n)\qquad\text{and}\qquad
\rho\left(\bigcup_{n=1}^\infty\Ker(f^n)\right) \subset 
           \bigcup_{n=1}^\infty\Ker(g^n)
\] 
This means that, in the decomposition $\cv=\cv'\oplus\cv''$ that
comes from $f$ and the decomposition $\cw=\cw'\oplus\cw''$ that
comes from $g$, we must have the compatibility that $\rho:\cv\la\cw$
must split as $\rho=\rho'\oplus\rho''$ for a unique choice of
$\rho':\cv'\la\cw'$ and $\rho'':\cv''\la\cw''$.
\ermk

\lem{L47.9}
Let $k$ be an algebraically closed field, let $Y$ is a projective curve over
$k$ with only simple nodes, and  
let $\cv^*$ be a cochain complex of vector bundle on $Y$. 
Then any cochain map $f^*:\cv^*\la\cv^*$
leads to a canonical decomposition $\cv^*=\cv^*_1\oplus\cv^*_2$
such that 
\be
\item
$f^*$ decomposes as $f^*=f_1^*\oplus f^*_2$ for maps $f^*_1:\cv^*_1\la\cv^*_1$
and $f^*_2:\cv^*_2\la\cv^*_2$.
\item
The map $f^*_1$ is an isomorphism while the map $f^*_2$ is 
locally nilpotent. By \emph{locally nilpotent} we mean that, for any 
integer $i\in\zz$, there exists an integer $n_i$ (which may depend on $i$)
for which $f^i_2:\cv^i_2\la\cv^i_2$ satisfies $(f^i_2)^{n_i}=0$.
\setcounter{enumiv}{\value{enumi}}
\ee
Moreover: if $f$ is null homotopic then the complex $\cv^*_1$ must be 
contractible.
\elem

\prf
The existence of the canonical decomposition is by Remark~\ref{R47.7};
in each degree $i$ the map $f^i:\cv^i\la\cv^i$ leads to a 
decomposition $\cv^i=\cv_1^i\oplus\cv^i_2$, and the differential
$\partial^i:\cv^i\la\cv^{i+1}$ must split as a direct
sum $\partial^i=\partial^i_1\oplus\partial^i_2$ for suitable
$\partial^i_1:\cv^i_1\la\cv^{i+1}_1$ and 
$\partial^i_2:\cv^i_2\la\cv^{i+1}_2$.

It remains to prove the ``moreover'' statement. Suppose therefore
that $f^*$ is null homotopic. Then the isomorphism 
$f_1^*:\cv_1^*\la\cv_1^*$ can be factored as the composite
\[\xymatrix{
\cv_1^* \ar[r]^-i & \cv^*\ar[r]^-f &\cv^*\ar[r]^-p &\cv_1^*
}\]
where $\cv_1^*\stackrel i\la\cv^*\stackrel p\la\cv^*_1$ are
the canonical inclusion and projection from the direct sum.
Thus the fact that $f$ is null homotopic forces the
isomorphism $f_1^*:\cv_1^*\la\cv_1^*$ to also be.
\eprf

\pro{P47.11}
Let $k$ be an algebraically closed field, let $Y$ be a projective curve over
$k$ with only simple nodes, and
let $\K\big[\vect Y\big]$ be the homotopy category of cochain
complexes of vector bundles
on $Y$. Then every idempotent in $\K\big[\vect Y\big]$
splits.
\epro

\prf
Choose an idempotent in $\K\big[\vect Y\big]$, and let it be represented
by the cochain map $e^*:\cv^*\la\cv^*$. The
assumption that $e^*$ is idempotent in $\K\big[\vect Y\big]$
means that
$e^*$ and $(e^*)^2$ are homotopic.

By Lemma~\ref{L47.9} we may decompose $\cv^*$, along the map
$e^*:\cv^*\la\cv^*$, as $\cv^*=\cv_1^*\oplus\cv^*_2$ is such a way that
\be
\item
The map $e^*$ may be written as $e_1^*\oplus e_2^*$, for cochain maps
$e_1^*:\cv_1^*\la\cv_1^*$ and $e_2^*:\cv_2^*\la\cv_2^*$.
\item
The map $e_1^*$ is an isomorphism
while the map $e_2^{*}$ is locally nilpotent.
\setcounter{enumiv}{\value{enumi}}
\ee
The fact that $(e^*)^2-e^*$ is null homotopic forces the 
direct summands $(e_1^*)^2-e_1^*$
and $(e_2^*)^2-e_2^*$ to be null homotopic.
Since $e_1^*$ is an isomorphism and $e_1^*(e_1^*-1)$ is null homotopic
we deduce that $e_1^*$ is homotopic to the identity.

It remains to show that $e_2^*$ is null homotopic. 
Replacing $e^*$ by $e^*_2$, we are reduced 
to showing 
\be
\setcounter{enumi}{\value{enumiv}}
\item
Suppose $e^*:\cv^*\la\cv^*$ is a cochain map, and assume that 
$e^*$ is locally nilpotent and that $e^*-(e^*)^2$ is null homotopic.
Then $e^*$ is null homotopic. 
\setcounter{enumiv}{\value{enumi}}
\ee
Observe the formal equalities
\begin{eqnarray*}
e^* &=& \sum_{i=1}^\infty \left[(e^*)^i-(e^*)^{i+1}\right]\\
    &=&\big[e^*-(e^*)^2\big]\cdot\left[\sum_{i=0}^\infty (e^*)^i\right]\\
\end{eqnarray*}
The infinite sums make sense since the local nilpotence guarantees
that in each degree the sums are finite. Thus we have produced a 
factorization of the cochain map $e^*:\cv^*\la\cv^*$ as the composite
of two cochain maps
\[\xymatrix@C+60pt{
\cv^*\ar[r]^-{\sum_{i=0}^\infty (e^*)^i} &\cv^* \ar[r]^-{e^*-(e^*)^2} & \cv^*
}\]
where the second is null homotopic. Hence $e^*$ is null homotopic.
\eprf

\rmk{R47.19}
The proof of Proposition~\ref{P47.11} generalizes, it 
works to show that 
idempotents split not only in the 
unbounded homotopy category $\K\big[\vect X\big]$, but 
also in the bounded subcategories 
$\K^b\big[\vect X\big]$, $\K^-\big[\vect X\big]$
and $\K^+\big[\vect X\big]$.
But for the categories $\K^b$, $\K^-$ and $\K^+$
the fact that idempotents split is 
not new: the reader can construct proofs using the 
methods of
\cite[Proposition~3.1]{Bokstedt-Neeman93}. See the
proof of \cite[Proposition~3.4]{Bokstedt-Neeman93} for 
an outline. For a later proof that is more {\it K--}theoretical
see Balmer and Schlichting~\cite[Section~2]{Balmer-Schlichting}.
More precisely: the idempotent-completeness of
$\K^-(\ce)=\D^-(\ce^\oplus)$ and of $\K^+(\ce)=\D^+(\ce^\oplus)$
may be found in \cite[Lemma~2.4]{Balmer-Schlichting},
while the idempotent-completeness of $\K^b(\ce)=\D^b(\ce^\oplus)$
follows from \cite[Theorem~2.8]{Balmer-Schlichting}.
\ermk

\section{The counterexample}
\label{S48}

Throughout this section $k$ will be a fixed algebraically closed
field. If $Y$ is a projective curve over $k$ with only simple nodes
as singularities, then $\vect Y$ will be the category of vector 
bundles over $Y$. And now we return to the {\it K--}theoretic implications
of what we have proved.

\rmk{R48.1}
The {\it K--}theoretic import of Proposition~\ref{P47.11} is that 
$\KK_{-1}\big[{\vect Y}^\oplus\big]=0$. This follows from
Schlichting~\cite[Corollary 6 of Section~9]{Schlichting06}. We
briefly remind the reader: in the category $\cm$
of model categories we have a commutative
square
\[\xymatrix@C+30pt{
\Ch^b\big[\vect Y\big]\ar[r]\ar[d] & \Ch^-\big[\vect Y\big]\ar[d] \\
\Ch^+\big[\vect Y\big]\ar[r]  & \Ch\big[\vect Y\big]
}\]
The functor $\D:\cm\la\ct$ takes this to the commutative
square
\[\xymatrix@C+30pt{
\K^b\big[\vect Y\big]\ar[r]\ar[d] & \K^-\big[\vect Y\big]\ar[d] \\
\K^+\big[\vect Y\big]\ar[r]  & \K\big[\vect Y\big]
}\]
and recalling that, by Reminder~\ref{R8.-1},
we have $\K^?(\ce)=\D^?(\ce^\oplus)$ for $?$ any of
$b$, $+$, $-$ or the empty restriction, this puts
us squarely in the situation of 
Schlichting~\cite[Corollary 6 of Section~9]{Schlichting06}.
The connective part of the
functor $\KK$ takes the commutative square
in $\cm$ to the homotopy cartesian square
\[\xymatrix@C+30pt{
\KK_{\geq0}\Big(\Ch^b\big[\vect Y\big]\Big)\ar[r]\ar[d] & 0\ar[d] \\
0\ar[r]  & \KK_{\geq0}\Big(\Ch\big[\vect Y\big]\Big)
}\] 
making $\KK_{\geq0}\Big(\Ch\big[\vect Y\big]\Big)$ a delooping of
$\KK_{\geq0}\Big(\Ch^b\big[\vect Y\big]\Big)
\cong\KK_{\geq0}\big[{\vect Y}^\oplus\big]$.
And $\KK_{-1}\big[{\vect Y}^\oplus\big]$ can therefore be computed
as $\KK_0$ of the idempotent completion
of $\D\Big(\Ch\big[\vect Y\big]\Big)=\K\big[\vect Y\big]$.
Proposition~\ref{P47.11} tells us 
that the triangulated category $\K\big[\vect Y\big]$ is idempotent-complete,
and since it has vanishing $\KK_0$ we deduce that
$\KK_{-1}\big[{\vect Y}^\oplus\big]=0$.
\ermk

\exm{E8.5}
As in Remark~\ref{R48.1} we let $Y$ be a projective curve over
$k$ with simple nodes as singularities.
Now recall the homotopy fibration of Proposition~\ref{P8.1}. With
$\ce=\vect Y$ we have an exact sequence 
$\KK_{-1}(\ce^\oplus)\la\KK_{-1}(\ce)\la\KK_{-2}(M')$, where
$\D(M')=\A^b(\ce)$ is the homotopy category of 
bounded acyclic cochain complexes.

Remark~\ref{R48.1} gives the vanishing of
$\KK_{-1}(\ce^\oplus)$, and  the sequence
of Proposition~\ref{P8.1} becomes  
$0\la\KK_{-1}(\ce)\la\KK_{-2}(M')$. This shows
that $\KK_{-2}(M')$ contains $\KK_{-1}(\ce)$ as a submodule.
But there are known examples of nodal curves with non-vanishing
$\KK_{-1}$; see Weibel~\cite[Exercise~III.4.4]{Weibel13}. Hence
there are examples of model categories $M'\in\cm$,
with non-vanishing negative \kth, and 
such that $\D(M')=\A^b(\ce)$ has a bounded \tstr.

Note: as presented in Weibel's book the nodal curves $U$
for which $\KK_{-1}\big[{\vect U}\big]\neq0$ are affine.
But the Mayer-Vietoris exact sequence, for assembling
$\KK\big[{\vect Y}\big]$ from
Zariski open covers, permits us
to pass to the compactification of these affine curves,
which can be chosen to have only simple
nodes as singularities and also have non-vanishing
$\KK_{-1}\big[{\vect Y}\big]$.
The Mayer-Vietoris sequence was proved to be exact 
by Weibel~\cite[Main Theorem]{Weibel89} for reduced 
quasiprojective varieties with isolated singularities, and 
in general by
Thomason and 
Trobaugh~\cite[Theorem~8.1]{ThomTro}.\footnote{For the discerning,
careful  reader
who checks the reference to \cite{ThomTro}:
the only schemes $Y$ considered
in the current article have been quasiprojective, they have ample
line bundles, this makes them satisfy the resolution
property,
and the natural map $\D^b\big[\vect Y\big]\la\dperf Y$
is an equivalence of categories.

Thomason and
Trobaugh~\cite{ThomTro} allows 
pathological schemes which do not satisfy the resolution property,
and for such schemes the right \kth\ to work with, the one
for which the Mayer-Vietoris sequence is exact,
is the \kth\ that comes from
perfect complexes and not the one corresponding
to $\Big(\Ch^b\big[\vect Y\big],\A^b\big[\vect Y\big]\Big)$. 
It isn't known 
how pathological a
scheme $Y$ has to be 
for this to be an issue. Non-separated schemes can fail
to satisfy the resolution property, but the reader 
can check Totaro~\cite[Question~1, page~3 of the Introduction]{Totaro04}: 
for all we know every scheme with an affine diagonal 
might have
the resolution property.}

The non-vanishing of  $\KK_{-1}\big[{\vect Y}\big]$
can also be deduced from the Mayer-Vietoris
sequence for the conductor square; see 
Pedrini and Weibel~\cite[Theorem~A.3]{Pedrini-Weibel94}.  
\eexm

\rmk{R48.9}
In Remark~\ref{R8.3} we promised the reader that, for the 
$\ce=\vect Y$ of Example~\ref{E8.5}, 
the category $\ce$ will be proved hereditary.
It is time to deliver on the promise.
 
If $\cv,\cv'$
are vector bundles on $Y$ then there is a spectral sequence,
converging to $\Ext^{p+q}(\cv,\cv')$, whose $E^2$ term has entries
$H^p\big[Y,\ce\mathit{xt}^q(\cv,\cv')\big]$. Now as $\cv$ is a locally
free sheaf the local Ext sheaves $\ce\mathit{xt}^q(\cv,\cv')$ vanish
when $q>0$. And $H^p\big[Y,\ce\mathit{xt}^0(\cv,\cv')\big]$
vanishes if $p>1$, just because $Y$ is a curve.
\ermk

\def\cprime{$'$}
\providecommand{\bysame}{\leavevmode\hbox to3em{\hrulefill}\thinspace}
\providecommand{\MR}{\relax\ifhmode\unskip\space\fi MR }
\providecommand{\MRhref}[2]{%
  \href{http://www.ams.org/mathscinet-getitem?mr=#1}{#2}
}
\providecommand{\href}[2]{#2}

\end{document}